# The Discovery of a New Series of Uniform Polyhedra.

Rinus Roelofs
Lansinkweg 28
7553AL Hengelo
The Netherlands
E-mail: rinus@rinusroelofs.nl
www.rinusroelofs.nl

**Abstract**

This article describes two new families of uniform polyhedra as well as the construction of the models of some of them. It is possible that these families has been discovered before, but I have been unable to find any publication about it. I welcome any information about earlier publication.

## 1. Introduction

**1.1. First discovery.** Can you imagine that the object in Figure 1 is a uniform polyhedron?

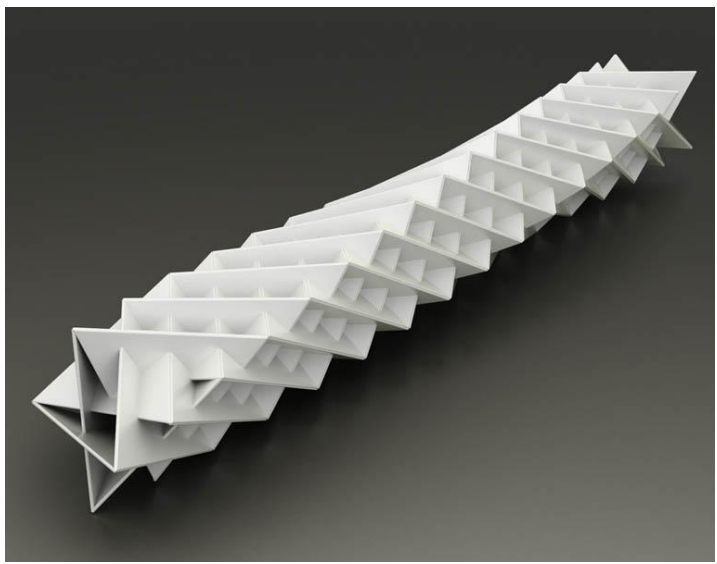

**Figure 1**: *Helical Star Deltahedron 9-4(2).*

After discovering a series of constructions of this kind, which I called *Helical Star Deltahedra* I had the idea that these objects are in fact uniform polyhedra. All these helical star deltahedra have two properties in common: they have only regular faces (all faces are equilateral triangles) and all vertices are congruent (at each vertex six triangles are joined together).

**1.2. Definition.** There are several different definitions in use to describe the family of regular polyhedra.
In the paper "An Enduring Error" by Branko Grünbaum [1] this matter is very clearly analysed. In the section about General Polyhedra he quotes Coxeter, Longuet-Higgins and Miller from [2]:
A *polyhedron* is a finite set of [planar] polygons such that every side of each belongs to just one other, with the restriction that no subset has the same property. The polygons and their sides are called faces and edges. The faces are not restricted to be convex, and may surround their centres more than once (as, for example, the pentagram, or five-sided star polygon, which surrounds its centre twice). Similarly, the faces at a vertex of a polyhedron may surround the vertex more than once. A polyhedron is said to be *uniform* if its faces are regular while its vertices are all alike. By this we mean that one vertex can be transformed into any other by a symmetry operation."

Polyhedra are surfaces composed of polygons such that each edge is adjacent to two polygons. Regular polyhedra are those that are composed of only one type of regular polygon. Because of the Coxeter's discovery of the infinite regular polyhedra clearly we can generalise the definition quoted above even more by removing the word "finite". After looking at my construction, with the definition in mind, I came to the conclusion that this might be an interesting discovery. To explain this step by step, we have to go back to history of the regular polyhedra.

## 2. Regular Polyhedra.

**2.1. Platonic Solids.** About 350 B.C. Plato decribed the five convex regular polyhedra (Figure 2) in his book *Timaeus*. It can easily be proven that there are exactly five convex regular polyhedra.

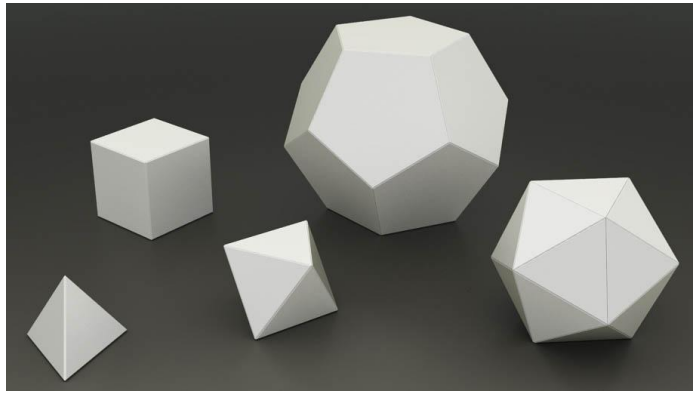

**Figure 2**: *The Platonic Solids.*

**2.2. Kepler-Poinsot.** In 1619 a first extension of the series of regular polyhedra was published by Kepler in his *Harmonices Mundi*. His two new regular polyhedra are not convex and are built with star-shaped faces Figure 3). About 200 years after Kepler's publication, Poinsot found two other non-convex regular polyhedra and these were added to the family. They are shown in Figure 4.

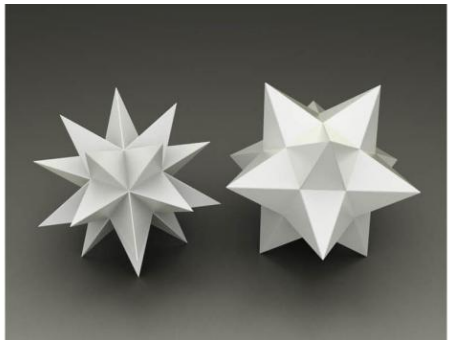
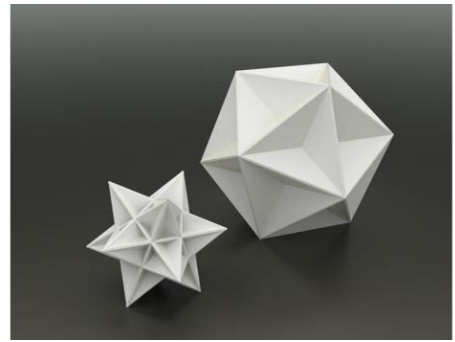

**Figure 3**: *Kepler polyhedra.*   **Figure 4**: *Poinsot polyhedra.*

**2.3. Coxeter.** In 1937 Coxeter, Longuet-Higgins and Miller described three new regular polyhedra to be added to the nine previously known [2]. These infinite regular polyhedra are shown in Figures 5-7.

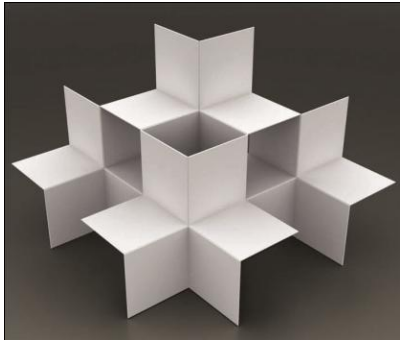 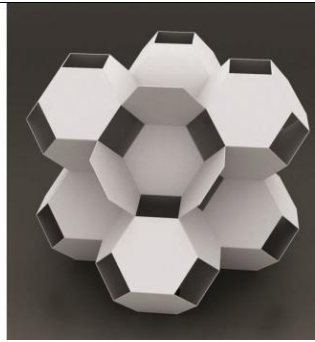 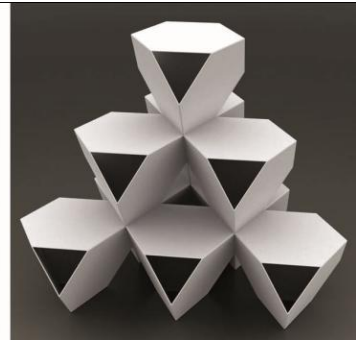

**Figure 5**: *Coxeter 4-6.*     **Figure 6**: *Coxeter 6-4 .*     **Figure 7**: *Coxeter 6-6.*

**2.4. Infinite Polyhedra.** In 1974 Wachman, Burt and Kleinmann published their book "Infinite Polyhedra" [4] that gave an introduction to and a description of the infinite uniform polyhedra. They divide the infinite polyhedra in three families: A. Cylindrical polyhedra, infinite in one direction only, B. Polyhedra with two noncollinear vectors of translation and C. Polyhedra with three noncoplanar vectors of translation. Here I will concentrate on family A.

### 3. Cylindrical Polyhedra

**3.1. Antiprisms.** Wachman, Burt and Kleinmann's book describes six infinite families of cylindrical polyhedra. If we do not allow two adjacent faces of a polyhedron to be coplanar then only two of their families are left: the family of polyhedra composed of antiprismatic rings and the family composed of helicoidal strips of equilateral triangles. In this publication the authors classify these polyhedra as regular uniform polyhedra.

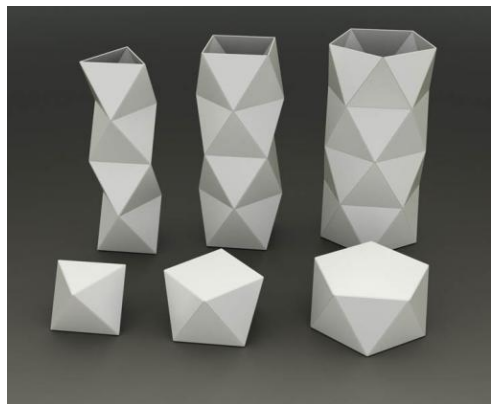

**Figure 8**: *Infinite polyhedra composed of antiprismatic rings.*

**3.2. Tetrahelix.** The tetrahelix can be seen as assemblage of tetrahedra and it is also the first member of the family of polyhedra composed of helicoidal strips. It can be unfolded as in Figure 11 in three side-by-side connected strips of equilateral triangles. In general, the polyhedra of this family, which we will call *helical deltahedra*, can be unfolded into a series of side-by-side connected strips.

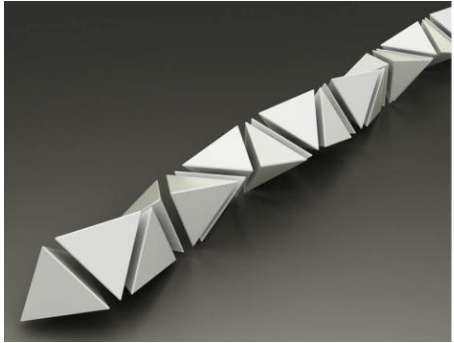 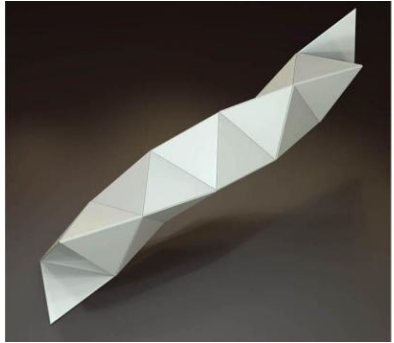 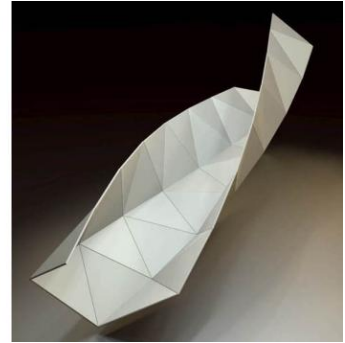

**Figure 9**: *Tetrahedra*   **Figure 10**: *Tetrahelix.*   **Figure 11**: *Unfolding the tetrahelix*

## 4. New Polyhedra.

**4.1. Strips of equilateral triangles.** To construct the helical deltahedra we can start with any number of side-by-side connected strips of equilateral triangles. In Figures 12-14 we show the process of folding a helical deltahedron with five strips. When we can connect the left edge of the left-most strip to the right edge of the right-most strip in such a way that the triangles become connected edge to edge, the final helical deltahedra is formed.

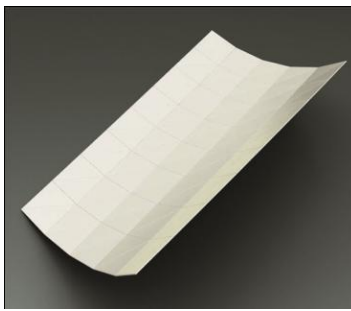 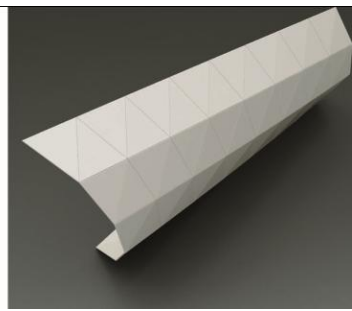 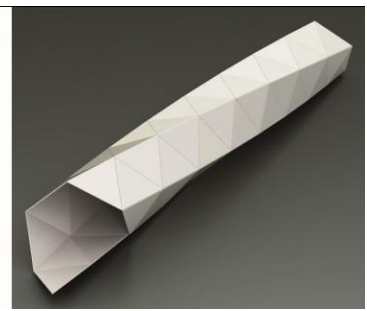

**Figure 12**: *Five strips.*   **Figure 13**: *Folding.*   **Figure 14**: *Helical deltahedron.*

After we have reached the situation in Figure 14, we can continue to make the folding edges sharper. The strips will now intersect as can be seen in Figure 15-16, and at a certain moment the original left-most edge can be connected to the original right-most edge again (Figure 17). The final new polyhedron is a uniform polyhedron. All faces are equal and the vertices are all congruent.

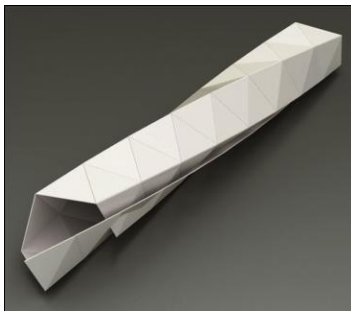 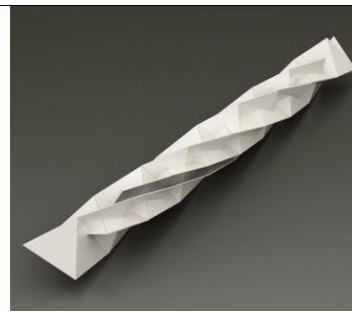 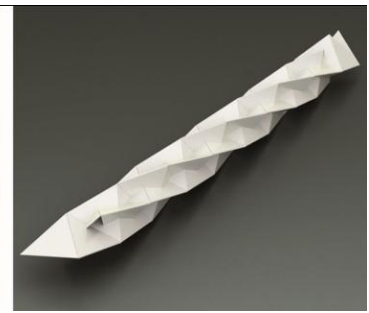

**Figure 15**: *Intersecting*   **Figure 16**: *Almost closed.*   **Figure 17**: *Helical star deltahedron.*

**4.2. Shift.** Depending on the number of strips you begin with, there is usually more than one way to connect the left-most edge of the set of strips to the right-most edge. There has to be a shift to get a helical deltahedron. But the number of steps in this shift may vary (Figure 18-20). Starting with five strips there are two possibilities to create a helical deltahedron in which faces do not intersect with other faces.

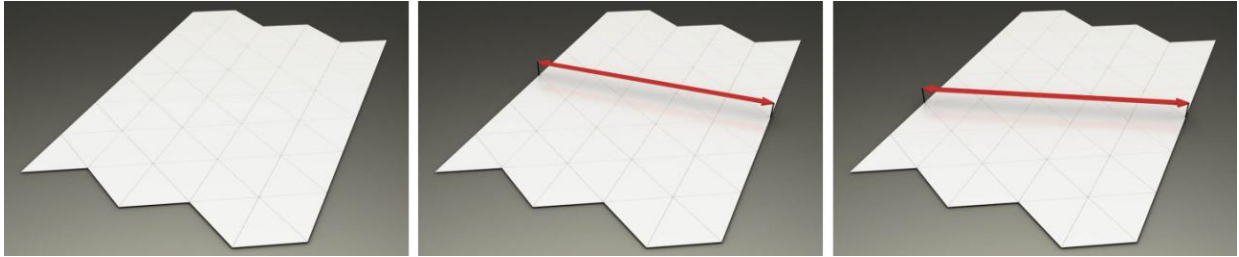

**Figure 18**: *5 strips of triangles*     **Figure 19**: *Connection a.*     **Figure 20**: *Connection b.*

Also for the new structure, with the intersection, there are two ways to create a regular polyhedron. They are shown in Figure 21 and 22. Because of the star shape in the polyhedron I will call this new family of uniform polyhedra "helical star deltahedra". Figure 23 shows an overview of regular star polygons is shown starting from 5-2 up to 12-5.

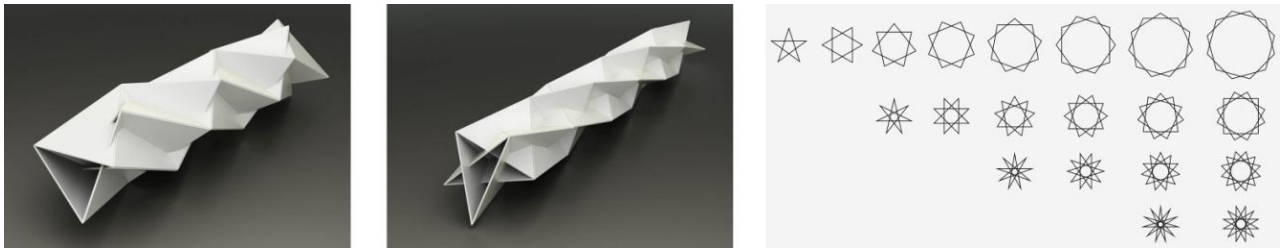

**Figure 21**: *Helical Star Deltahedron 5-2(1).*   **Figure 22**: *Helical Star Deltahedron 5-2(2)*     **Figure 23**: *Stars.*

**4.3. Enumeration.** The series of regular helical star deltahedra is infinite because we can start with any number of strips of equilateral triangles. Depending on the number of strips one or more star shapes are possible. And also the possible number of steps in the shifts depends on the number of strips. Starting with five strips we can make two different regular helical star deltahedra: each has a star shape, combined with a one or a two-step shift (Figures 24-25).

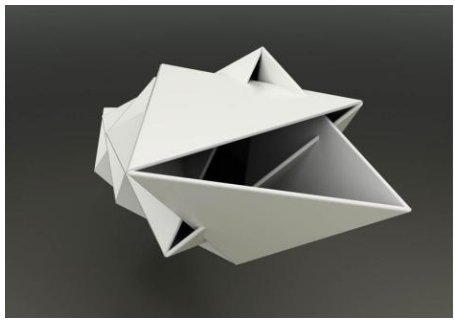 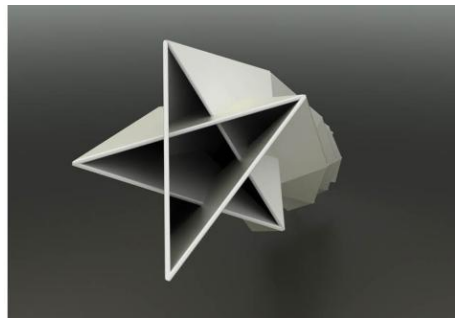

**Figure 24**: *Star 5-2, one step.*     **Figure 25**: *Star 5-2, two steps .*

In naming the regular helical star deltahedra we give the star shape (this also gives the number of strips) followed by the number of steps in the shift. Six strips yields only one regular helical star deltahedron. The structure we get when we apply a shift of two steps is a compound of two tetrahelices. I have found 64 different regular helical star deltahedra for stars up to the 12 point stars.

## 5. Examples

**5.1. Infinite series.** The number of regular helical star deltahedra is infinite. From five to twelve strips however we count 64. Four examples are shown in Figures 26-29. I was surprised by the variety and the beauty of these polyhedra and therefore a have started making real models of these shapes.

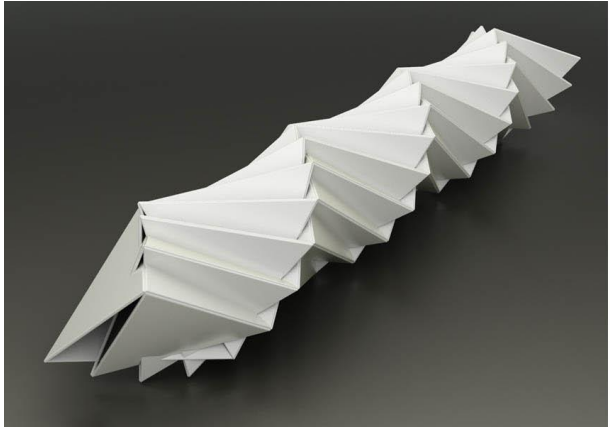
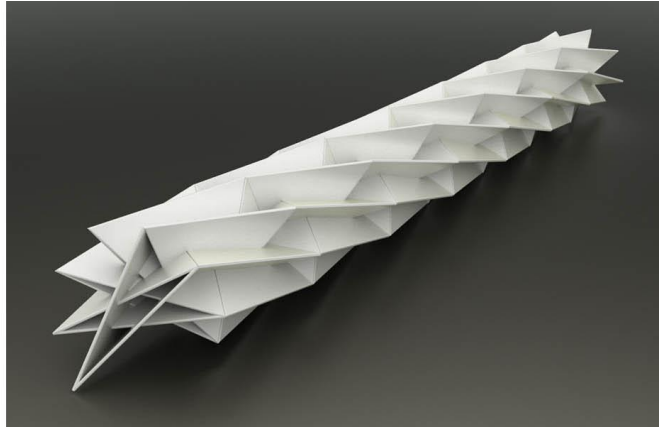

**Figure 26**: *Regular Helical Star Deltahedron 9-4(1.)*   **Figure 27**: *Regular Helical Star Deltahedron 9-4(4).*

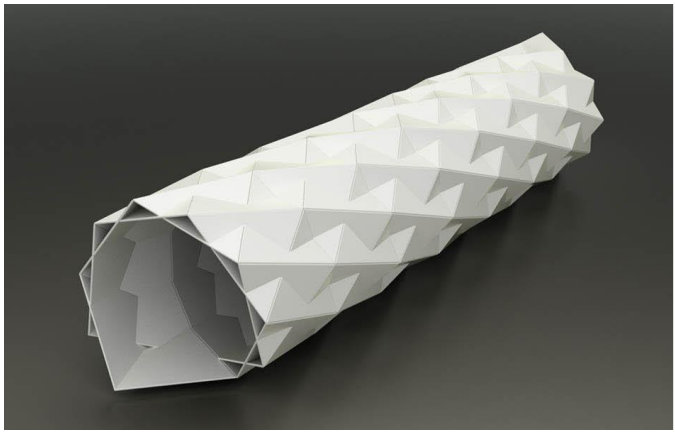
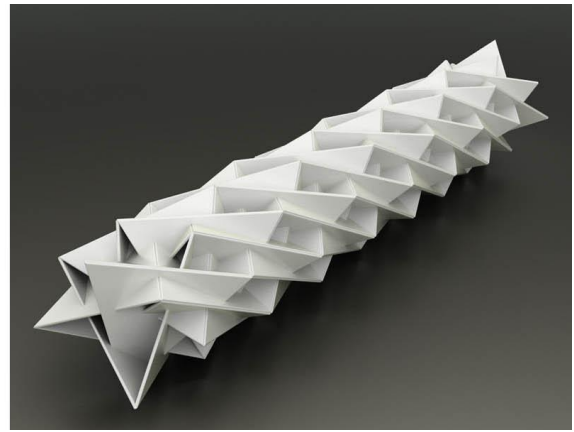

**Figure 28**: *Regular Helical Star Deltahedron 11-2(2).*   **Figure 29**: *Regular Helical Star Deltahedron 12-5(3).*

**5.2. Edges.** To get an impression of the internal structure of the polyhedra I made a number of models with only the edges of the equilateral triangles. Figures 30-33 show examples. The model in Figure 31 is constructed with one type of connector (all the vertices are congruent) and pipes of one length only.

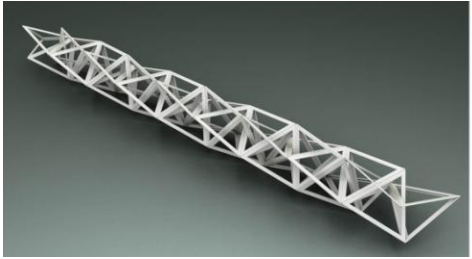 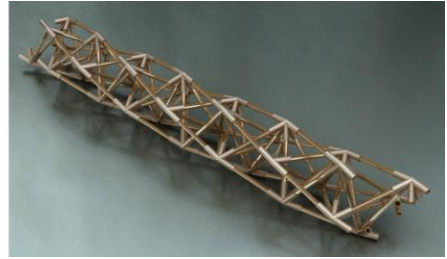

**Figure 30**: *Edges of 5-2(2).*     **Figure 31**: *Frame of 5-2(2).*

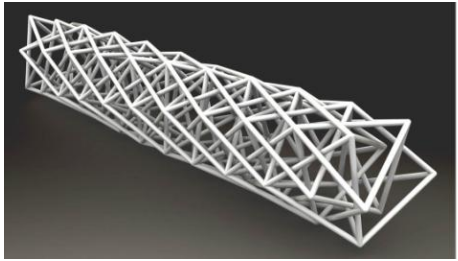 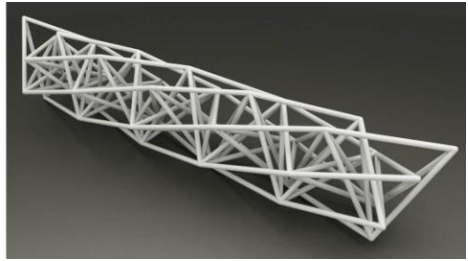

**Figure 32**: *Edges of 8-3(1).*     **Figure 33**: *Edges of 7-2(2).*

## 6. Alternative Construction

**6.1. Sliding.** After making a number of models, I saw another possible way of constructing models of some of the regular helical star polyhedra. With a laser cutter I made paper modules of two connected equilateral triangles that could be slid together. This turned out to work fine for the 5-2(2) (regular helical star polyhedra based on the 5-2 star, shift 2 steps) (Figures 34-36) as well as for the 7-2(2) (Figure 37).

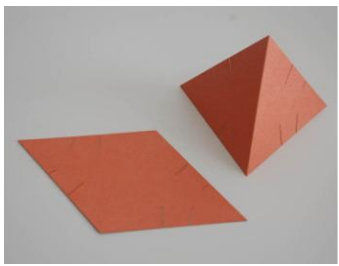 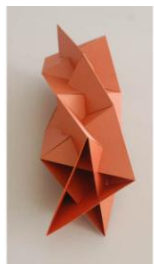 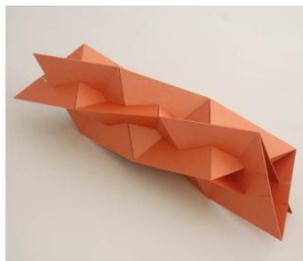 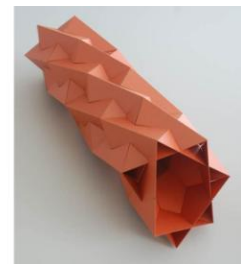

 **Figure 34**: *Modules.*     **Figure 35**: *5-2(2).*    **Figure 36**: *5-2(2).*     **Figure 37**: *7-2(2).*

**6.2. Stars 8-3 and 7-3.** Even more complex members of the family can be made using this method. Figures 38-40 show the model of the regular helical star deltahedron 8-3(3) and Figures 41-42 show the 7-3(3).

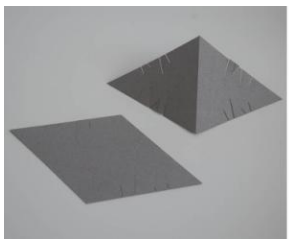 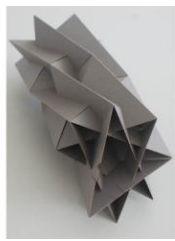 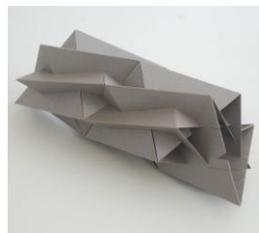 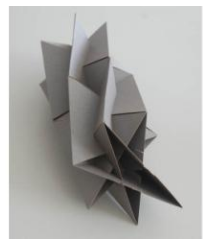 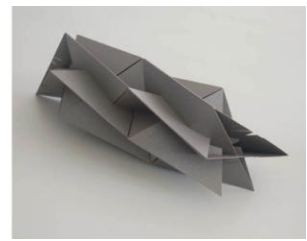

**Figure 38**: *Modules.* **Figure 39**: *8-3(3).*  **Figure 40**: *8-3(3).*     **Figure 41**: *7-3(3).*     **Figure 42**: *7-3(3).*

# 7. Vertex Figures

**7.1. Vertex Figure.** When we look at the vertex figure of the regular helical star deltahedra shown so far we see that it is a non-convex hexagon without intersections (Figure 43a).

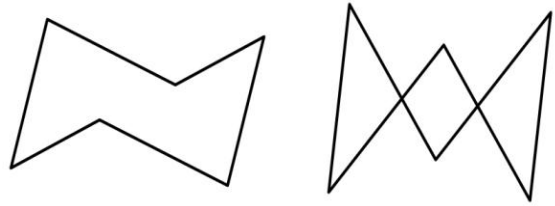

**Figure 43a**: *Vertex Figure without intersection*     **Figure 43b**: *Vertex Figure with intersection*

After seeing this I investigated the possibility of constructing regular helical star deltahedra which have vertex figures with intersections as in Figure 44b. It turned out to be possible and so far I have constructed a first set of 12 different regular helical star deltahedra with vertex figures of the type of Figure 43b.

**7.2. Examples.** Two examples of this second family are shown in Figures 44 and 45. Classification of members of this family is more complicated because we cannot construct these polyhedra with the same folding technique as in section 4.1. This will be one of the subjects for further research.
The computer program Rhinoceros and especially the plug-in Grasshopper, developed by David Rutten, were a great help in this research project.

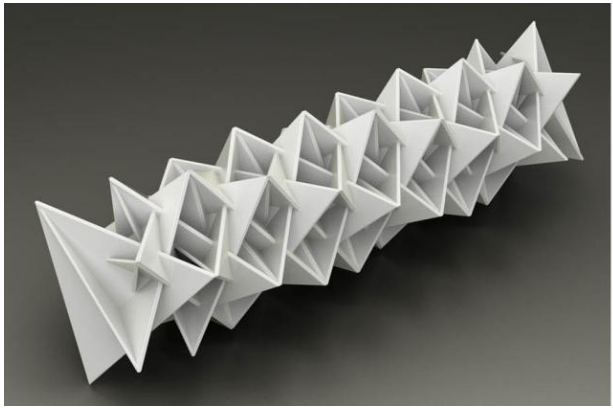
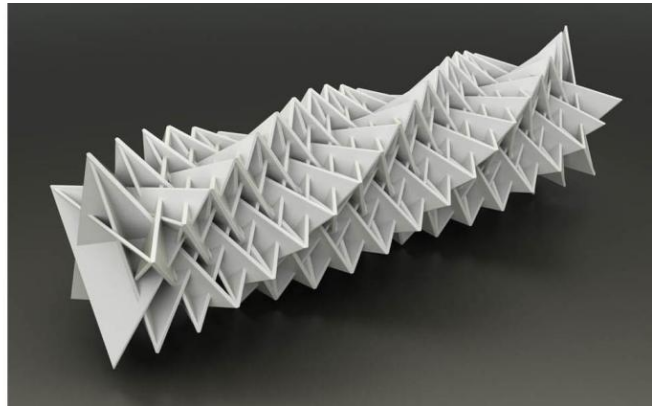

**Figure 44**: *Example 1.*     **Figure 45**: *Example 2.*

**7.3. Further Research.** A more detailed version of the paper is being prepared for publication in an appropriate mathematical journal.